\documentclass{proc-l}

\usepackage{latexsym,amssymb,times,mathptm,enumerate,amsthm}
\usepackage{amsthm}

\usepackage{hyperref}

\usepackage{mathbbol}

\def\proof{\noindent{\bf{Proof.} }}
\def\sqr#1#2{{\vcenter{\hrule height.#2pt
        \hbox{\vrule width.#2pt height#1pt \kern#1pt
                \vrule width.#2pt}
        \hrule height.#2pt}}}

\def\tratto{\mbox{\rule{2mm}{.2mm}$\;\!$}}

\newtheorem{theorem}{Theorem}[section]

\newtheorem{corollary}[theorem]{Corollary}
\newtheorem{lemma}[theorem]{Lemma}
\newtheorem{proposition}[theorem]{Proposition}

\newtheorem{remark}[theorem]{Remark}

\newtheorem{example}[theorem]{Example}

\newcommand{\m}{\mathfrak{m}}

\newcommand{\ol}[1]{\ensuremath{\overline{#1}}}

\numberwithin{equation}{section}

\newcommand{\ds}{\displaystyle}
\DeclareMathOperator{\Spec}{\rm{Spec}}
\DeclareMathOperator{\Max}{\rm{Max}}
\newcommand{\rees}[1]{\ensuremath{\mathcal{R}{(#1)}}}

\begin{document}

\title[Reductions of ideals]{Generators of reductions of ideals in a local Noetherian ring with finite residue field}


\author[L. Fouli]{Louiza Fouli}
\address{Department of Mathematical Sciences,  New Mexico State University,  Las Cruces, New Mexico 88003, USA}
\email{lfouli@nmsu.edu}
\thanks{The first author was partially supported by a grant from the Simons Foundation, grant \#244930.}

\author[B. Olberding]{Bruce Olberding}
\address{Department of Mathematical Sciences,  New Mexico State University,  Las Cruces, New Mexico 88003, USA}
\email{olberdin@nmsu.edu}




\commby{Irena Peeva}

\begin{abstract}
Let $(R,\m)$ be a local Noetherian ring with residue field $k$. 
While much is known about the generating sets of reductions of ideals of $R$ if $k$ is infinite, the case in which $k$ is finite is less well understood. We investigate the existence (or lack thereof) of proper reductions of an ideal of $R$ and the number of generators needed for a reduction in the case $k$ is a finite field. When $R$ is one-dimensional, we give a formula for the smallest integer $n$ for which every ideal has an $n$-generated reduction. It follows that in a one-dimensional local Noetherian ring every ideal has a principal reduction if and only if the number of maximal ideals in the normalization of the reduced quotient of $R$ is  at most $|k|$. In higher dimensions, we show that for any positive integer, there exists an  ideal of $R$ that does not have an $n$-generated reduction and that if $n \geq \dim R$ this ideal can be chosen to be $\m$-primary. In the case where $R$ is a two-dimensional regular local ring, we construct an example of an integrally closed $\m$-primary ideal that does not have a $2$-generated reduction and thus answer in the negative a question raised by Heinzer and Shannon.

\end{abstract}

\maketitle


\section{Introduction}

Let $R$ be a (commutative) Noetherian ring, and let $I$ be an ideal of $R$. A {\it reduction} of $I$ is a subideal $J$ of $I$ such that $I^{n+1} = JI^n$ for some $n>0$; equivalently, $J$ is a subideal of $I$ such that $\overline{I} = \overline{J}$, where $^{^{\tratto{}}}$ denotes the integral closure of the corresponding ideal in $R$.  Northcott and Rees \cite{NR} proved that if $R$ is a local Noetherian ring with infinite residue field and Krull dimension $d$, then every ideal of $R$ has a $d$-generated reduction, that is, a reduction that can be generated by $d$ elements. This result and its generalizations involving analytic spread underlie many of the applications of the theory of reductions to local algebra. For example, reductions and the analytic spread have been instrumental in describing the asymptotic properties of an ideal $I$ of $R$, the Cohen-Macaulay property of the Rees algebra $\rees{I}$ of $I$, and the  blowup Proj $\rees{I}$ of $\Spec{R}$ along the subscheme defined by $I$. However, if the residue field of $R$ is finite, then there may exist ideals of $R$ that do not have a $d$-generated reduction, and so the applicability of reductions in the case of local rings with finite residue field is more limited. 

In this article, we examine the extent to which the result of Northcott and Rees involving $d$-generated reductions fails  in the case of finite residue field. 
We prove two main results, the first of which is  devoted to one-dimensional rings and the second to rings of higher dimension. In the one-dimensional case we find the optimal choice for  replacing $d$ in the result of  Northcott and Rees with the smallest possible positive integer, an integer that depends only on the size of the residue field and the cardinality $|\Max(\overline{R_{red}})|$ of the set of maximal ideals in the normalization $\overline{R_{red}}$ of the reduced quotient $R_{red} = R/\sqrt{0}$ of $R$. 
 
\smallskip

{\noindent}{\bf Theorem A}.   {\it Let $R$ be a one-dimensional local Noetherian ring with finite residue field $k$. 
The smallest positive integer $n$ for which every ideal of $R$ has an $n$-generated reduction is 
$$ n = \lceil \: -1 +  \log_{|k|}\left(\: |k| + (|k|-1)\cdot |\Max(\overline{R_{red}})|\: \right) \: \rceil.$$
 }

For a one-dimensional local Noetherian ring $R$ of multiplicity $e$, every ideal of $R$ can be generated by $e$ elements \cite[Theorem 1.1, p.~49]{Sa}. Thus, the number $n$ in the theorem is at most $e$.

If $R$ is a one-dimensional local Noetherian ring with infinite residue field, then every ideal of $R$ has a principal reduction. Using Theorem A, we extend this result to  one-dimensional local rings with residue field of any size (see Corollary~\ref{count}).

\medskip

{\noindent}{\bf Corollary}.   {\it Let $R$ be a one-dimensional local Noetherian ring with residue field $k$. 
Every ideal of $R$ has a principal reduction if and only if $|\Max(\overline{R_{red}})| \leq |k|$.
 }  

\medskip

In particular, if $R$ is a complete local Noetherian domain, then every ideal of $R$ has a principal reduction (Corollary~\ref{complete case}).  Thus  in dimension one there are interesting  local Noetherian rings {\it with  finite residue field}  for which 
 the result of Northcott and Rees holds, i.e., every ideal of $R$ has a reduction generated by $\dim(R)$-many elements.  Moreover, in dimension one, even if there are ideals without a principal reduction, we are at least guaranteed the existence of a bound on the number of elements needed to generate a reduction.  

Moving beyond dimension one, we use Theorem A to  show (see Theorem~\ref{general case}) that unlike in the case of infinite residue field, 
 no such bound exists for a local Noetherian ring with  finite residue field and  dimension at least $2$:

\medskip

{\noindent}{\bf Theorem B.} 
{\it    Let $(R, \m)$ be a local  Noetherian ring of  dimension $d \geq 2$. If  the residue field of $R$ is finite, then for each positive integer $n$ there is an ideal of $R$ that is minimally generated by $n$ elements and does not have a proper reduction. If also $n\geq d$, this ideal can be chosen to be $\m$-primary. 
} 

\medskip

{\it Notation.} Throughout the article, $Q(R)$ denotes the total quotient ring of the ring $R$, $\overline{R}$ is the integral closure of $R$ in $Q(R)$, and $\overline{I}$ is the integral closure of the  ideal $I$. We denote by $R_{red}$  the reduced ring $R/\sqrt{0}$, where $\sqrt{0}$ is the nilradical of $R$. The set of maximal ideals of $R$ is denoted $\Max(R)$. 
\section{Preliminaries} 

In this section we  develop a criterion for when  every ideal in a local Noetherian ring with finite residue field has an  $n$-generated reduction. This criterion, Proposition~\ref{survive}, will be used in the proofs of the main results in the next section.  
The first  lemma, which is a routine application of well-known properties of reductions, concerns the transfer of  reductions in a ring $R$ to its reduced quotient $R_{red}$.

\begin{lemma} \label{n lemma}
Let $R$ be a ring and let $n$ be a positive integer. A finitely generated ideal $I$  of $R$ has an $n$-generated reduction if and only if $IR_{red}$ has an $n$-generated reduction. 
\end{lemma}

\proof 
It is clear that if every ideal of $R$ has $n$-generated reduction, then every  ideal of $R_{red}$ has an $n$-generated reduction. Conversely, let $I$ be an ideal of $R$ and suppose that $J$ is an $n$-generated ideal of $R$ such that $JR_{red}$ is a reduction of $IR_{red}$. Then $J \subseteq I + \sqrt{0}$.  Write $J = (x_1,\ldots,x_n)R$.  For each $i$, there is $y_i \in I$ and $z_i \in \sqrt{0}$ such that $x_i = y_i + z_i$.  Thus $J + \sqrt{0} = (y_1,\ldots,y_n)R + \sqrt{0}$.  Let $K = (y_1,\ldots,y_n)R$.  Then $K \subseteq I$ and $KR_{red} = JR_{red}$.  Thus $KR_{red}$ is a reduction of $IR_{red}$. We claim that $K$ is a reduction of $I$. 
Since 
$\overline{KR_{red}} = \overline{IR_{red}}$ 
we have by \cite[Proposition~1.1.5]{HS} that  $\overline{K}R_{red} = \overline{I}R_{red}$. By \cite[Remark~1.1.3(5)]{HS}, the nilradical $\sqrt{0}$ of $R$ is contained in every integrally closed ideal of $R$, so we conclude that $\overline{K} = \overline{I}$. Thus $K$ is an $n$-generated reduction of $I$. 
\qed

\medskip

As discussed in the introduction, if $(R, \mathfrak{m})$ is a one-dimensional local Noetherian ring with infinite residue field, then  every ideal of $R$ has a principal reduction. Removing the restriction to infinite residue field, we can assert in general that every ideal has a principal reduction if and only if every $\mathfrak{m}$-primary ideal has a principal reduction, or more generally:

\begin{proposition} \label{height 0} Let $(R,\mathfrak{m})$ be a one-dimensional local Noetherian ring, and let $n$ be a positive integer.  If every $\mathfrak{m}$-primary ideal of $R$ has an $n$-generated reduction, then every ideal of $R$  has an $n$-generated reduction. 

\end{proposition} 

\proof  Let $I$ be an ideal of $R$. 
  By Lemma~\ref{n lemma}, it suffices to show that $IR_{red}$ has an $n$-generated reduction, so we may assume without loss of generality that $R$ is a reduced ring and $I$ is a proper ideal of $R$. Let $A = (0:I)$.   Since $R$ is local, reduced and one-dimensional, the  ideal $I + A$ is $\mathfrak{m}$-primary. Indeed, suppose that $P$ is a minimal prime of $R$ with $I + (0:I) \subseteq P$. Then $R_P$ is a field as $R$ is reduced and hence $IR_P \subseteq PR_P=0$. Then $R_P=(0R_P : IR_P) \subseteq PR_P$, a contradiction. Since $I+A$ is $\mathfrak{m}$-primary, then by assumption $I+A$ has an $n$-generated reduction $J$, say $(I+A)^{k+1} = J(I+A)^k$ for some $k>0$.   Using the fact that $IA = 0$, we have  $ I^{k+1} + A^{k+1} = (I+A)^{k+1} = J(I+A)^k = J(I^k + A^k) = JI^k + JA^{k}$.
  Write $J = (x_1+a_1,\ldots,x_n+a_n)R$, where each $x_i \in I$ and each $a_i \in A$. 
  Then $I^{k+1} + A = JI^k + A = (x_1,\ldots,x_n)I^k + A$.  
 Since $(I \cap A)^2 =0$ and $R$ is reduced, we have $I \cap A = 0$.  Therefore, from $ I^{k+1} + A = (x_1,\ldots,x_n)I^k + A$, we conclude that $I^{k+1}  = (x_1,\ldots,x_n)I^k$, which proves that $I$ has an $n$-generated reduction. 
\qed

\medskip

The next lemma and proposition 
give criteria for when every ideal in a Noetherian ring $R$ has an $n$-generated reduction. The stronger result, Proposition~\ref{survive}, requires that $R$ is also reduced, local and one-dimensional. In light of Theorem B, the one-dimensional assumption is necessary in the proposition.  To state  Lemma~\ref{n+1}, we recall that the arithmetic rank, ara$(I)$, of a proper ideal $I$ of a Noetherian ring $R$ is the least number $n$ such that $\sqrt{I} = \sqrt{(x_1,\ldots,x_n)R}$ for some $x_1,\ldots,x_n \in I$.

\begin{lemma} \label{n+1} Let $R$ be a Noetherian ring, let $I$ be a proper ideal of $R$ and let $n$ be an integer such that $n \geq $ {\rm ara}$(I)$.   If  each $(n+1)$-generated ideal $J \subseteq I$ with $\sqrt{J} = \sqrt{I}$ has an $n$-generated reduction, then $I$ has an $n$-generated reduction.

\end{lemma}

\proof
Let $I$ be a proper ideal of $R$ and let  $${\mathcal F} = \{\overline{J}:  J \subseteq I, \ \sqrt{J} = \sqrt{I} \mbox{ and }J \mbox{ is an } n \mbox{-generated ideal of } R\}.$$  
Since $R$ is a Noetherian ring and $n \geq $ ara$(I)$, the set ${\mathcal F}$ is nonempty and  contains a  maximal element $\overline{J}$, where $J$ is an $n$-generated ideal with $J \subseteq I$ and $\sqrt{J} = \sqrt{I}$. Suppose that $\overline{J} \subsetneq \overline{I}$. Then $I \not \subseteq \overline{J}$ and we may choose $y \in I \setminus \overline{J}$. Now ${J+yR}$ is an $(n+1)$-generated ideal, and hence by assumption there is an $n$-generated ideal $K$ of $R$ such that $K \subseteq J+yR$ and $\overline{K} = \overline{J+yR}$.  But then $\overline{J} \subseteq \overline{K} \in {\mathcal F}$ and the maximality of $\overline{J}$ in ${\mathcal F}$ forces $\overline{J} = \overline{J+yR}$, a contradiction to the fact that $y \not \in \overline{J}$. Therefore, $\overline{J} = \overline{I}$ and  $J$ is an $n$-generated reduction of $I$.   
\qed

\begin{proposition} \label{survive} Let $R$ be a reduced one-dimensional local Noetherian ring and let $n$ be a positive integer. Then the following are equivalent:
\begin{enumerate}[$(1)$]
\item 
Every    ideal of  $R$ has an $n$-generated  reduction.

\item For all  $x_1,\ldots,x_{n+1} \in \overline{R}$ for which $(x_1,\ldots,x_{n+1})\overline{R} = \overline{R},$ there is an $n$-generated $R$-submodule  of $(x_1,\ldots,x_{n+1})R$ that does not survive in $\overline{R}$. 
\end{enumerate}
\end{proposition}

\proof
In the proof we use the fact  that every ideal of $\overline{R}$ is a principal ideal of $\overline{R}$.   This can be seen as follows. 
   By \cite[Theorem 12.3, p.~66]{Huc}, the fact that $R$ is a reduced local  Noetherian ring of dimension $\leq 2$ implies $\overline{R}$ is a finite product of Noetherian integrally closed domains. Since also $\dim \overline{R} = \dim R = 1$,  we have that $\overline{R}$ is a finite product of Dedekind domains. The fact that $\overline{R}$ is semilocal implies these Dedekind domains are principal ideal domains. As    a finite product of principal ideal domains, the ring $\overline{R}$ has the property that every ideal is a principal ideal.

 First,
suppose every  ideal of $R$ has an $n$-generated reduction and let $x_1,\ldots,x_{n+1} \in \overline{R}$ such that $(x_1,\ldots,x_{n+1})\overline{R} = \overline{R}$. Let $A = (x_1,\ldots,x_{n+1})R$. 
 Then there is a nonzero divisor $r \in R$ such that $I:=rA$ is an ideal of $R$.    
Since $A\overline{R} = \overline{R}$, we have $AQ(R) = Q(R)$, and hence $A$ contains a nonzero divisor of $R$.  Therefore, $I=rA$ also contains a nonzero divisor of $R$.
  By assumption, there is
an $n$-generated reduction $J$ of $I$.  
 Let $m>0$ be such that  $I^{m+1} = JI^{m}$.   
 As we have established, every ideal of $\overline{R}$ is a principal ideal. Thus $I\overline{R}$ is a principal ideal of $\overline{R}$ that is 
 necessarily generated by a nonzero divisor of $\overline{R}$ since $I$ contains a nonzero divisor. 
  Since $I^{m+1}\overline{R} = JI^m \overline{R}$ and principal ideals generated by a   nonzero divisor admit no proper reductions, we obtain 
 $I\overline{R} = J\overline{R}$.  Therefore, $\overline{R} = A\overline{R} = r^{-1}I\overline{R}= r^{-1}J\overline{R}$.  Since $J \subseteq  I$, we have $r^{-1}J \subseteq A$, and hence $r^{-1}J$ is an $n$-generated $R$-submodule of $A$ that does not survive in $\overline{R}$.

Conversely, suppose that for all  $x_1,\ldots,x_{n+1} \in \overline{R}$ with $(x_1,\ldots,x_{n+1})\overline{R} = \overline{R},$ there is an $n$-generated $R$-submodule of $(x_1,\ldots,x_{n+1})R$  that does not survive in $\overline{R}$.  
  To prove that every  ideal of $R$ has an $n$-generated reduction, it suffices by Lemma~\ref{n+1} to show that every $(n+1)$-generated  ideal of $R$ has an $n$-generated reduction, since ${\rm{ara}} (I)=1$ for any proper non-zero ideal $I$ of $R$. Let $a_1,\ldots,a_{n+1} \in R$ and let $I = (a_1,\ldots,a_{n+1})R$. By Proposition~\ref{height 0}, it suffices to consider the case in which $I$ is $\m$-primary, where $\m$ is the maximal ideal of $R$.

    Since every ideal of $\overline{R}$ is a principal ideal,  we have $I \overline{R} = t\overline{R}$ for some  $t \in \overline{R}$. 
    Since $I$ is $\m$-primary and $R$ is reduced, $I$ contains a nonzero divisor, so $t$ is a nonzero divisor in $\overline{R}$.  
     For each $i \in \{1,\ldots,n+1\}$, let $x_i = a_it^{-1}$.  Then 
 $\overline{R} = t^{-1}t \overline{R} = t^{-1}I \overline{R} = (x_1,\ldots,x_{n+1})\overline{R}$. 
 By
 assumption  there exists an $n$-generated $R$-submodule $A$ of $(x_1,\ldots,x_{n+1})R$ that does not survive in $\overline{R}$.  
 Since we have $(x_1,\ldots,x_{n+1})\overline{R} = \overline{R} = A\overline{R}$, it follows that  
 $$I\overline{R} = t(x_1,\ldots,x_{n+1})\overline{R} = tA \overline{R}.$$ Moreover, $tA \subseteq  t(x_1,\ldots,x_{n+1})R = I$.  By \cite[Proposition~1.6.1]{HS} the fact that $I\overline{R} = tA\overline{R}$ implies
  $\overline{I} = \overline{tA}$. Since $tA \subseteq I$, 
 this proves that $tA$ is an $n$-generated reduction of $I$. \qed

\section{Main results}\label{Main results}

In this section we prove the main results of the paper. After proving the first theorem, which deals with the one-dimensional case, we indicate how Theorem A of the introduction follows. At the end of the section in Theorem~\ref{general case}, we prove Theorem B of the introduction. 

\begin{theorem} \label{no red} Let $R$ be a one-dimensional local Noetherian ring with finite residue field $k$, and let $n$ be a positive integer. Then every ideal of $R$ has an $n$-generated reduction if and only if 
$$|\Max(\overline{R_{red}})| \leq \frac{|k|^{n+1} -|k|}{|k|-1}.$$   
\end{theorem}

\proof  
 By Lemma~\ref{n lemma}, every ideal of $R$ has an $n$-generated reduction if and only if every ideal of $R_{red}$ has an $n$-generated reduction. 
Thus  
 it suffices to prove the theorem in the case where $R$ is a reduced ring. 
Throughout the proof, we let $U$ denote a set of $|k|$-many elements of $R$ such that $R/\m = \{u + \m:u \in U\}$, where $\m$ is the maximal ideal of $R$. 
   We assume $0 \in U$. Since $R$ is a local ring, all the nonzero elements of $U$ are units in $R$.  
     We denote the elements of the Cartesian  product $U^n$ by ${\bf u} = (u_1,\ldots,u_n)$. 
     Let $${\mathcal J} =   \{(i,{\bf u}) \in \{1,\ldots,{n+1}\} \times U^n: u_j = 0 {\mbox{ for all }} j \geq i \}.$$
     Then 
     $$     |{\mathcal J}|  =   1 + |k| + |k|^2 + \cdots + |k|^{n} 
      =  \frac{|k|^{n+1} - 1}{|k|-1}.$$

Now we prove the theorem. 
Suppose first that $$|\Max(\overline{R})| > \frac{|k|^{n+1} -|k|}{|k|-1}.$$   We show there is an ideal of $R$ that does not have an $n$-generated reduction. 
 By Proposition~\ref{survive}, it suffices to show that there are $x_1,\ldots,x_{n+1} \in \overline{R}$ such that $\overline{R} = (x_1,\ldots,x_{n+1})\overline{R}$ and every $n$-generated  $R$-submodule of $(x_1,\ldots,x_{n+1})R$ survives in $\overline{R}$.

         By assumption, $$|\Max(\overline{R})| \geq   \frac{|k|^{n+1} -|k|}{|k|-1} + 1 = \frac{|k|^{n+1} - |k| + |k| - 1}{|k|-1} = |{\mathcal J}|.$$ 
   Therefore, we may index a set $\{M_{i,{\bf u}}: (i,{\bf u}) \in {\mathcal J}  \} $ of $|{\mathcal J}|$-many maximal ideals of $\overline{R}$ by ${\mathcal J}$.  Since the ideals $M_{i,{\bf u}}$ of $\overline{R}$ are maximal, the diagonal map $$\phi: \overline{R}  \longrightarrow \prod_{(i,{\bf u}) \in {\mathcal J}} \: \overline{R}/M_{i,{\bf u}}:x \mapsto (x+M_{i,{\bf u}})$$
is a surjective ring homomorphism. 
 Thus we may choose  $x_1,\ldots,x_{n+1} \in \overline{R}$ such that for each $(i,{\bf u}) \in {\mathcal{J}}$ and $j = 1,\ldots,n+1$, we have 
$$ x_j + M_{i,{\bf u}} 
\: = \: \begin{cases} 
  0 + M_{i,{\bf u}} & \textrm{ if \: $i < j$}, \\
     1 + M_{i,{\bf u}} & \textrm{ if \: $i=j$}, \\
       u_j + M_{i,{\bf u}}  & \textrm{ if \: $i > j$}. \\
   \end{cases} $$

Let $K = (x_1,\ldots,x_{n+1})R$.  We claim that $\overline{R} = (x_1,\ldots,x_{n+1})\overline{R}$ and every $n$-generated  $R$-submodule of $K$ survives in $\overline{R}$. To show this we first prove that $K/\m K$ has dimension $n+1$ as a $k$-vector space.  Indeed, suppose $i\in \{1,\ldots, n+1\}$ and   $$x_i + \m K = \sum_{j \ne i}r_j x_j  + \m K$$ for some $r_j \in R$.  For each $j \neq i$, the choice of  $x_j$ implies  $x_j \in M_{i,{\bf 0}}$ (with ${\bf 0} = (0,\ldots,0)$). 
Since $$x_i - \sum_{j \ne i}r_j x_j  \in \m K \subseteq M_{i,{\bf 0}},$$ we conclude that $x_i \in M_{i,{\bf 0}}$, contrary to the fact that by the choice of $x_i$, we have $x_i + M_{i,{\bf 0}} = 1 + M_{i,{\bf 0}}$.
This contradiction shows that $K/\m K$ has dimension $n+1$ as a $k$-vector space.   In particular, $K$ cannot be generated as an $R$-module by fewer than $n+1$ elements.

Let $I$ be an $n$-generated $R$-submodule of $K= (x_1,\ldots,x_{n+1})R$. We claim that $I$ survives in $\overline{R}$.  
 Since $I$ can be generated as an $R$-module by $n$ elements, the dimension of the $k$-vector space $I/\m I$ is at most $n$.  Since $K$ cannot be generated by fewer than $n+1$ elements, by adding as many of the elements $x_i$ to $I$ as needed we can assume without loss of generality  that $I$ is an $R$-submodule of $K$ such that $I$ can be generated by $n$ but no fewer  elements. In particular,
 the $k$-dimension of $I/\m I$ is $n$. 
 Then Nakayama's Lemma implies that $I$ is  generated by $n$ elements $y_1,\ldots,y_n$ of the form 
 \begin{center}
 $y_j = u_{j,1}x_1+u_{j,2}x_2 + \cdots + u_{j,n+1}x_{n+1},$ \:
where $u_{j,i} \in U$.  
\end{center}
Since $I/\m I$ is a $k$-vector space of  dimension $n$, the $n \times (n+1)$ matrix $$(u_{j,i}+\m:j=1,\ldots,n, \: i = 1,\ldots,n+1),$$ whose entries are in the field $k$, has rank $n$. Elementary row operations produce a rank $n$ matrix in reduced row echelon form such that for some $i =1,\ldots,n~+~1$, deleting the $i$-th column yields the $n \times n$ identity matrix. It follows from this observation and Nakayama's Lemma that there are $u_1,\ldots,u_{i-1} \in U$ such that 
 $$I  = (\{x_j - u_jx_i: j < i\} \cup \{x_{i+1},\ldots,x_{n+1}\})R.$$  

Let ${\bf u} = (u_1,\ldots,u_{i-1},0,\ldots,0) \in U^n$ so that $(i,{\bf u}) \in {\mathcal J}$.
We claim that $I \subseteq M_{i,{\bf u}}$.  First, observe that $x_{i+1},\ldots,x_{n+1} \in M_{i,{\bf u}}$ by the choice of these elements, so  
it remains to show that for each $j < i$, we have  $x_j - u_jx_i \in M_{i,{\bf u}}$. Let $j < i$ and notice that by the choice of $x_i$, we have
\begin{center} $x_i + M_{i,{\bf u}} = 1 + M_{i,{\bf u}}$ \:\:  and  \:\: 
 $x_j + M_{i,{\bf u}} = u_j + M_{i,{\bf u}}$.\end{center}
  Therefore, $$x_j  - u_jx_i + M_{i,{\bf u}} = u_j - u_j\cdot 1 + M_{i,{\bf u}} = 0 + M_{i,{\bf u}},$$ so that $x_j-  u_jx_i \in M_{i,{\bf u}}$, proving the claim that $I \subseteq M_{i,{\bf u}}$.  We conclude that $I$ survives in $\overline{R}$. This shows that every $n$-generated $R$-submodule of $K = (x_1,\ldots,x_{n+1})R$ survives in $\overline{R}$. Therefore, $R$ has an ideal that does not have an $n$-generated reduction. 
 
 Conversely, 
suppose that 
$$|\Max(\overline{R})| \leq \frac{|k|^{n+1} -|k|}{|k|-1}.$$   To prove that every ideal of $R$ has an $n$-generated reduction, it suffices by Proposition~\ref{survive} to show that for all $x_1,\ldots,x_{n+1} \in \overline{R}$ with $(x_1,\ldots,x_{n+1})\overline{R} = \overline{R}$, there  is an $n$-generated $R$-submodule $I$ of $(x_1,\ldots,x_{n+1})R$ such that $I\overline{R} = \overline{R}$. 

 Let $x_1,\ldots,x_{n+1} \in \overline{R}$ such that $(x_1,\ldots,x_{n+1})\overline{R} = \overline{R}$.
 For each $(i, {\bf u}) \in {\mathcal J}$ with ${\bf{u}} = (u_1,\ldots,u_n)$, consider the $n$-generated $R$-submodule of $(x_1,\ldots,x_{n+1})R$ defined by  
 $$I_{i,{\bf u}} =  (x_1-u_1x_i, \: x_2 - u_2x_i,\ldots, x_{i-1}-u_{i-1}x_i, \: x_{i+1},\ldots,x_{n+1})R.$$
 
 We claim first that if $(s,{\bf u}), (t,{\bf v}) \in {\mathcal J}$ and $I_{s,{\bf u}}$ and $I_{t,{\bf v}}$ are contained in a common maximal ideal of $\overline{R}$, then $s = t$ and ${\bf u} = {\bf v}$.  Suppose $M$ is a maximal ideal of $\overline{R}$ with $I_{s,{\bf u}} + I_{t,{\bf v}}\subseteq M$. 
 
 Suppose by way of contradiction that  $s \ne t$.  Without loss of generality, we can assume that $s < t$.  Then $x_t \in I_{s,{\bf u}}$, 
 so that since $$x_1-v_1x_t, x_2 - v_2x_t,\ldots, x_{t-1}-v_{t-1}x_t \in I_{t,{\bf v}} \subseteq M,$$ we have $x_1,x_2, \ldots,x_{t-1} \in M$.  Also, since $s < t$, we have $x_{t},\ldots,x_{n+1} \in I_{s,{\bf v}} \subseteq M$.  But then $\overline{R} = (x_1,\ldots,x_{n+1})\overline{R} \subseteq M$, a contradiction that implies $s = t$.  
 
 Next we claim that ${\bf u} = {\bf v}$.  Since $s = t$, we have $$x_1-u_1x_t,\ldots, x_{t-1}-u_{t-1}x_t \in I_{s,{\bf u}} \subseteq M.$$
 Similarly,  we have 
$$x_1-v_1x_t,\ldots, x_{t-1}-v_{t-1}x_t \in I_{t,{\bf v}} \subseteq M.$$
Therefore, for each $i < t$, we have $$x_i-u_ix_t - (x_i -v_ix_t) = (v_i-u_i)x_t \in M.$$
 An argument similar to the one in the preceding paragraph shows that $x_t \not \in M$, since otherwise every $x_j \in M$, for  $j=1,2,\ldots,n+1$, a contradiction. Thus, for each $i < t$ we have  $v_i - u_i \in M \cap R = \m$, where the last equality follows from the fact that as a maximal ideal of an integral extension of $R$, $M$ lies over $\m$.  Consequently,  for each $i < t$ we have $u_i +\m = v_i + \m$, which, since $u_i,v_i \in U$, forces $u_i = v_i$.  Since this holds for all $i< t$, we conclude ${\bf u} = {\bf v}$.  This proves that if $I_{s,{\bf u}}$ and $I_{t,{\bf v}}$ are contained in a common maximal ideal of $\overline{R}$, then $s = t$ and ${\bf u} = {\bf v}$. 
 
 Next, since no two distinct $R$-submodules of $(x_1,\ldots,x_{n+1})R$ from the set $\{I_{s,{\bf u}}:(s,{\bf u}) \in {\mathcal J}\}$ are contained in the same maximal ideal of $\overline{R}$, it follows that either one of these $R$-submodules does not survive in $\overline{R}$ or there are at least $|{\mathcal J}|$-many maximal ideals of $\overline{R}$.  As established at the beginning of the proof,  $$|{\mathcal J}|  
      =  \frac{|k|^{n+1} - 1}{|k|-1},$$ so if every $I_{s,{\bf u}}$  survives in $\overline{R}$, we conclude that  $$  \frac{|k|^{n+1} - 1}{|k|-1} \: \leq \: |\Max(\overline{R})| \: \leq \: \frac{|k|^{n+1} - |k|}{|k|-1},$$ 
      where the last inequality is given by assumption. This 
implies $$|k|^{n+1} - 1 \: \leq \: |k|^{n+1} - |k|,$$ which is impossible since $|k|>1$.  
 This contradiction implies some $I_{s,{\bf u}}$  does not survive in $\overline{R}$. In particular, some $n$-generated $R$-submodule of $(x_1,\ldots,x_{n+1})R$ does not survive in $\overline{R}$. Therefore, by Proposition~\ref{survive} every ideal of $R$ has an $n$-generated reduction. 
\qed 

\smallskip

Theorem A now follows easily from  Theorem~\ref{no red}: In the setting of the theorem, we seek the smallest positive integer $n$ such that 
$$|\Max(\overline{R_{red}})| \: \leq \: \frac{|k|^{n+1} -|k|}{|k|-1};$$   that is, we need the smallest positive integer $n$ such that  
$$|k| + (|k|-1) \cdot |\Max(\overline{R_{red}})| \: \leq \: {|k|^{n+1}}.$$   Equivalently, 
$$\log_{|k|}\left(|k| + (|k|-1) \cdot |\Max(\overline{R_{red}})|\right) \: \leq \: {n+1},$$
which yields the conclusion of Theorem A.

\begin{remark} \label{min max excellent} {\em If $R$ is a one-dimensional Noetherian local  domain, then 
 the number of  maximal ideals of $\overline{R}$  is the same as the number  of  minimal prime ideals in $\widehat{R}$ \cite[Corollary 5]{Kat}. Moreover, if $R$ is a reduced local Noetherian ring with geometrically regular formal fibers, then the number of maximal ideals of $\overline{R}$ is the same as the number of minimal prime ideals of the completion $\widehat{R}$ of $R$ \cite[Theorem 6.5]{Die}. Thus for excellent local Noetherian rings $R$, the bound in Theorem~\ref{no red}  can be restated using the minimal primes of the completion of $R_{red}$ rather than the maximal ideals of the normalization of $R_{red}$.  }
\end{remark}

Next we give a criterion for the existence of principal reductions in the one-dimensional case.

\begin{corollary} \label{count} Let $R$ be a one-dimensional local Noetherian ring with residue field~$k$. Then every  ideal of $R$ has a principal reduction if and only if $$|\Max(\overline{R_{red}})| \leq |k|.$$
\end{corollary}

\proof Apply Theorem~\ref{no red} in the case $n=1$.    
\qed

\begin{corollary} \label{mult} Let $R$ be a one-dimensional local Noetherian  domain with residue field $k$. If  the  multiplicity of $R$ is at most $|k|$, then every ideal of $R$ has a principal reduction.
\end{corollary} 

\proof
 The completion of a  one-dimensional local Noetherian domain of multiplicity $e$ has at most $e$ minimal prime ideals (e.g., this follows from the multiplicity formula given in \cite[Theorem 14.7]{Mat}). 
 By Remark~\ref{min max excellent}, $|\Max(\overline{R})| \leq e \leq |k|$.  
 Thus the corollary is a consequence of  Corollary~\ref{count}.\qed 

\begin{remark} \label{Bass} {\em It follows from Corollary~\ref{mult} that if a one-dimensional Noetherian local domain  $R$ has multiplicity 2, then every ideal of $R$ has a principal reduction. This is known already for other reasons. In a one-dimensional local Cohen-Macaulay ring of  multiplicity $e$,   every ideal can be generated by $e$ elements \cite[Theorem 1.1, p.~49]{Sa}. Since $R$ has multiplicity $2$,  every ideal of $R$ is $2$-generated. 
In \cite[Theorem 3.4]{SV}, Sally and Vasconcelos prove that a ring in which every ideal is $2$-generated has the property that every $\m$-primary ideal has a principal reduction of reduction number at most $1$; see also \cite[Lemma~1.11]{Li}.}
\end{remark}

\begin{corollary}  \label{two minimal} Let $(R,\m)$ be  a one-dimensional Noetherian local domain whose completion has at most two minimal prime ideals. Then every  ideal of $R$ has a principal reduction.
\end{corollary}

\proof
Apply Remark~\ref{min max excellent} and  Corollary~\ref{count}.
\qed

\begin{corollary} \label{complete case} Every ideal of a one-dimensional complete local domain has a principal reduction.
\end{corollary}

\proof Apply   Corollary~\ref{two minimal}. \qed

\medskip

Next we use the one-dimensional case in Theorem~\ref{no red}  to show the absence of a bound on the number of generators of reductions in higher dimensions.

 \begin{theorem} \label{general case}
 Let $(R, \m)$ be a local  Noetherian ring of  dimension $d \geq 2$. If  the residue field of $R$ is finite, then for each positive integer $n$ there is an ideal of $R$ that is minimally generated by $n$ elements and does not have a proper reduction. If also $n\geq d$, this ideal can be chosen to be $\m$-primary.

 \end{theorem}

\proof Let $n > 0$ and let $k$ denote the residue field of $R$. For the first assertion in the theorem, it suffices to show that there exists an ideal of $R$ that is generated by $n$ elements and  has no $(n-1)$-generated reduction, see for example \cite[Proposition~8.3.3]{HS}.
Since $R$ is Noetherian and $R$ has   dimension $d>1$, there are infinitely many prime ideals of $R$ of dimension one. Choose a positive integer $t$ with 
$$t \: > \:
 \frac{|k|^{n} -|k|}{|k|-1}.$$
Let $P_1,\ldots,P_t$ be distinct dimension one prime ideals of $R$ and let $$A = R/(P_1 \cap \cdots \cap P_t).$$  Then $A$ is  a reduced one-dimensional local Noetherian ring with $t$ minimal prime ideals. As such, $\overline{A}$ is a direct product of $t$  integrally closed domains \cite[p.~64]{Huc}. Consequently, $\overline{A}$ 
has at least $t$ maximal ideals.  
By Lemma~\ref{n+1} and Theorem~\ref{no red}   there is an $n$-generated ideal $I$ of $A$ that does not have an $(n-1)$-generated reduction. 
 Let $J$ be an $n$-generated   ideal of $R$
  such that  
     $JA = I$. 
     If $J$ has an $(n-1)$-generated reduction $K \subseteq J$, then there is $m>0$ such that  $J^{m+1} = J^mK$. But then $I^{m+1} = J^{m+1}A = J^mKA = I^mKA$, so that $KA$ is an $(n-1)$-generated reduction of $I$, contrary to the choice of $I$. We conclude that the $n$-generated ideal $J$ has no $(n-1)$-generated  reduction. 
     Hence  $J$ is minimally generated by $n$ elements and has no proper reduction.

         It remains to prove the second assertion in the theorem. Suppose $n \geq d$.  
         If $n = d$, then  we may choose  any $n$-generated $\m$-primary ideal of $R$. By Krull's height theorem such an ideal must exist and cannot have an $(n-1)$-generated reduction since such a reduction would be an $\m$-primary ideal generated by $d-1$ elements.  
         Thus, if $n = d$, the proof is complete. 
         
         Assume  $n>d$.  
        With $I$ as above, choose $J$ to be an ${\m}$-primary ideal of $R$ with 
                   $P_1 \cap \cdots \cap P_t \subseteq J$ and $JA = I$.  
     Since $J$ is ${\m}$-primary, Krull's height theorem implies that the arithmetic rank of $J$ is $d$. Therefore, since $n -1 \geq d$, Lemma~\ref{n+1} implies that there is an $n$-generated ideal $J' \subseteq J$ such that $\sqrt{J'} = {\m}$ and  $J'$ does not have an $(n-1)$-generated reduction. 
         This shows that $J'$ is an $n$-generated $\m$-primary ideal with no $(n-1)$-generated reduction, so we conclude that $J'$ is minimally generated by $n$ elements and $J'$ has no  proper reduction. 
\qed

\begin{remark}\label{Question HeS} {\em  In \cite[Example 2.3]{HRR} an example is given of a two-dimensional Cohen-Macaulay local ring $R$ such that $R$ has finite residue field and the maximal ideal $\m$ of $R$ fails to have a 2-generated reduction. In \cite[Example 5.4]{HeS}, Abhyankar showed that certain canonically defined $\m$-primary ideals of a two-dimensional regular local ring $(R,\m)$ have $2$-generated reductions. 
This motivates the question of Heinzer and Shannon \cite[Question 5.6]{HeS} of whether every integrally closed $\m$-primary ideal of  a two-dimensional regular local ring $(R,\m)$ with finite residue field has a 2-generated reduction. Theorem~\ref{general case} guarantees the existence of $\m$-primary ideals without 2-generated reductions, but the ideals produced in the proof need not be integrally closed. In the next section we give an example that answers this question in the negative; see Example~\ref{m-primary int closed no red}. 

}
\end{remark}

\section{Examples}~\label{Examples}

In this section we give several examples to illustrate some of the ideas in Section~\ref{Main results}. In the first example, we show that in order for each ideal of 
 a one-dimensional local Noetherian domain to have a principal reduction, it is not sufficient that every integrally closed ideal has a principal reduction.

\begin{example}\label{counter} {For each  integer $n\geq 1$, there exists a one-dimensional local Noetherian domain $R$ such that every integrally closed ideal has a principal reduction, yet $R$ has an ideal that is minimally generated by $n$ elements and does not have a proper reduction.}
{\em Let $n \geq 1$. If $n =1$, then any one-dimensional local Noetherian domain suffices as an example, since a nonzero principal ideal does not have a proper reduction. Suppose $n>1$.  
 Choose $2^{n} - 1$ maximal ideals $M_1,\ldots,M_{2^{n}-1}$ of the ring ${\mathbb{Z}}_2[x]$, and let $S =  {\mathbb{Z}}_2[x] \setminus (M_1 \cup \cdots \cup M_{2^{n}-1})$.   Let $J$ denote the Jacobson radical of the localization ${\mathbb{Z}}_2[x]_S$ of ${\mathbb{Z}}_2[x]$ at the multiplicatively closed set $S$, and 
let $R = {\mathbb{Z}}_2 + J$.  Then $R$ is a one-dimensional local Noetherian domain with maximal ideal $J$ and  normalization  $\overline{R} =  {\mathbb{Z}}_2[x]_S$; see for example \cite[Lemma~1.1.4 and Proposition~1.1.7]{FHP}. Thus $\overline{R}$ is  a  PID with $|\Max(\overline{R})| = 2^{n}-1$.  Observe that $(J:_{Q(R)} J) = \overline{R}$, and thus since $\overline{R}$ is a PID, the  blow up of the maximal ideal $J$ of $R$ in the sense of Lipman \cite[p.~651]{Li} is $\overline{R}$. 
This ring has  the property that each localization at a maximal ideal has its embedding dimension (which is 1) equal to its multiplicity. Therefore, by \cite[Theorem~2.2]{Li}, 
 every integrally closed ideal of $R$ has a principal reduction of reduction number at most $1$. 
However, since the residue field $k$ of $R$ has two elements and $$|\Max(\overline{R})| = |k|^{n} - 1 > \frac{|k|^{n} -|k|}{|k|-1},$$ Lemma~\ref{n+1} and Theorem~\ref{no red} imply that $R$ has an ideal that is minimally generated by $n$ elements and does not have an $(n-1)$-generated reduction. By \cite[Proposition 8.3.3]{HS}, this ideal does not  have a proper reduction. 
}
\end{example}

Corollary~\ref{mult} implies that the maximal ideal of a one-dimensional local Noetherian domain of multiplicity $2$ has a principal reduction. 
The next example, which appears in \cite[Example~8.3.2]{HS}, shows that there exists a one-dimensional local Cohen-Macaulay ring of multiplicity $3$ whose {\it maximal ideal} does not have a principal reduction.

\begin{example} \label{HS example} {\rm (See \cite[Example~8.3.2]{HS})
Let $R$ be the ring given by $$R =  \ds\frac{{\mathbb{Z}}_2[[x,y]]}{(xy(x+y))}.$$ The multiplicity of $R$ is $3$. Let $\m$ denote the maximal ideal of $R$, and let $x', y'$ denote the images of $x$ and $y$ in $R$, respectively.  Notice that ${\rm Min}(R)=\{(x'),(y'),(x'+y')\}$ and by Remark~\ref{min max excellent} and Corollary~\ref{count} we are guaranteed that there exists an ideal that does not have any principal reductions. We claim that $\m$ is such an ideal. Indeed, suppose that $\m$ has a principal reduction $J=(f)$. Write $f = r_1x' + r_2y'$, for some $r_1,r_2 \in R$. 
Since $R = {\mathbb{Z}}_2 + \m$, there are $a_1,a_2 \in {\mathbb{Z}}_2$ such that for each $i$,  $r_{i} - a_{i} \in \m$. Let $f' = a_1x'  + a_2y'$, and let $J' = (f')$.  Then   $J \subseteq J' +\m I \subseteq I$ and thus $\ol{J'}=\ol{I}$ by  \cite[Lemma~8.1.8]{HS}. In this way, 
 we may assume that any principal reduction of $\m$ is generated by an element of the form $f=a_1x'+a_2y'$, with $a_1, a_2 \in {\mathbb{Z}}_2$. Thus,  the only possible generators for principal reductions of $\m$ are $x',y',x'+y'$, but each of these is in a minimal prime ideal of $R$. However, any  
reduction of $\m$ will be an $\m$-primary ideal, and thus $\m$ has no principal reductions.

}
\end{example}

Although the ring $R$ in Example~\ref{HS example} is not a domain, it can be used to produce similar examples that are  domains. 
We recall 
 a theorem of Lech \cite[Theorem 1]{Lech}: 
A complete local Noetherian ring $R$ with maximal ideal $\m$ 
is the completion of a local Noetherian domain if and only if 
(i) $\m =0$ or $\m \not \in$ Ass$(R)$,    
 and 
(ii) no nonzero integer of $R$ is a zero divisor. 

\begin{example} A one-dimensional local Noetherian domain whose maximal ideal does not have a principal reduction. {\rm  Let $R$ be as in Example~\ref{HS example}. Since $R$ is Cohen-Macaulay, $\m \not \in $Ass$(R)$, and since the only nonzero integer of $R$ is $1$, it follows from the theorem of Lech that there exists a local Noetherian domain $A$ with completion $R$.  As in the preceding example, the maximal ideal $\m$ of $R$ does not have a principal reduction. 
Since the maximal ideal of $R$ is extended from the maximal ideal of $A$, it follows that  the maximal ideal of $A$ does not have a principal reduction. Moreover, since $R$ is a complete intersection of multiplicity $3$, so is $A$.
}
\end{example}

\smallskip

The following example was suggested to us by Bill Heinzer.

\begin{example}\label{Bill}{\rm
Let $R=\mathbb{Z}_2[[x,y,z]]/(z^2-xy,xy(x+y)(x+y+z))$. One can show that ${\rm{Min}}(R)=\{(x,z),(y,z), (x+z,y+z),(x+y+z)\}$. Every 
linear form belongs to some minimal prime of $R$ and hence there is no principal reduction of the maximal ideal $\m=(x,y,z)$. On the other hand, the number of minimal primes is $4$ and thus the number of maximal ideals of $\overline{R_{red}}$ is $4$ by Remark~\ref{min max excellent}. 
%
Since also  $|k| <4$,   Corollary~\ref{count} and Remark~\ref{min max excellent} imply that there exists an ideal of $R$ that does not have a principal reduction. In this case, $\m=(x,y,z)R$ is such an ideal.}
\end{example}

The next example was suggested  by Bernd Ulrich.
\begin{example}\label{Bernd}
{\rm Let $R=\mathbb{Z}_2[[x,y,z]]/(z^2-xy,xy(x+y))$. One can show that ${\rm{Min}}(R)=\{(x,z),(y,z), (x+z,y+z)\}$. It is straightforward to show that $J_1=(x+y+z)$ is the only principal reduction of $\m=(x,y,z)$. The number of minimal primes is $3$ and even though $|k| <3$ it is the case that  $\m$ has a principal reduction. By Remark~\ref{min max excellent} and Corollary~\ref{count} there must exist a different ideal of $R$ that does not have a principal reduction. One can verify that $J_2=(x,y)$ is such an ideal. Moreover,  $J_2$ is another proper reduction of $\m$ and there is no principal reduction contained in $J_2$. Hence there are  minimal reductions of $\m$ having minimal generating sets of different sizes.}

\end{example}

All the examples in this section have been devoted to the one-dimensional case. In the next example, we illustrate the failure of $2$-generated reductions for the maximal ideal of a two-dimensional ring. 
As mentioned in Remark~\ref{Question HeS}, Heinzer, Ratliff and Rush  \cite[Example 2.3]{HRR} have given for each finite field $F$ an example   
of a two-dimensional Cohen-Macaulay local ring $R$ such that $R$ has residue field $F$, the associated graded ring of $\m$ is Cohen-Macaulay, and the maximal ideal $\m$ of $R$ fails to have a 2-generated reduction. For the case $F = {\mathbb{Z}}_2$,  we give a simple example of this phenomenon   that has the same  properties.

\begin{example}\label{2-dim}{\rm
Let $R=\mathbb{Z}_2[[x,y,z]]/(xy(x+y)(x+y+z))$. Then  $R$ is a two-dimensional Cohen-Macaulay ring and  the associated graded ring of $\m$, ${\rm gr}_{\m}(R)=\bigoplus \limits_{i=0}^{\infty}\m^{i}/\m^{i+1}$, is Cohen-Macaulay. We claim that $\m=(x,y,z)R$ does not have a 2-generated reduction. 
Indeed, as in Example~\ref{HS example}, if $\m$ has a 2-generated reduction, then one can assume that the two generators of the reduction are linear forms. Since the residue field is $\mathbb{Z}_2$ there are only seven linear forms in $x,y,$ and $z$. It is straightforward to check that any ideal in $R$ generated by any two  linear forms has height 1 and therefore cannot be a reduction of $\m$.
}
\end{example}

In Remark~\ref{Question HeS} we mentioned a question raised by Heinzer and Shannon. In \cite{HeS} they ask if every integrally closed, $\m$-primary ideal of a two-dimensional regular ring $(R, \m)$ with finite residue field has a $2$-generated reduction. While Theorem~\ref{general case} guarantees the existence of $\m$-primary ideals without $2$-generated reductions, it asserts nothing about whether there are such examples that are integrally closed. However, following  the proof of Theorem~\ref{general case}, we can use similar ideas to produce such an example.
To this end, we first construct in Example~\ref{m-primary no red} an explicit example of a 3-generated ideal $I$ in a two-dimensional regular local ring such that $I$ has no 2-generated reductions and  is not integrally closed.  

\begin{example}\label{m-primary no red}
{\rm Let $R=\mathbb{Z}_2[x,y]_{(x,y)}$ be a regular local ring with $\m=(x,y)$  the unique maximal ideal of $R$. The algorithm in the proof of Theorem~\ref{general case} produces the following $\m$-primary ideal of $R$ that does not have any $2$-generated reductions:
 $$I=(x^2y+xy^2, \: xy^5+xy^4+xy^3+x^3, \: y^8+xy^3+x^3+xy^2).$$ This is an $\m$-primary ideal. If there is a $2$-generated reduction $J$ of $I$, then we may assume that $J$ is generated by two ${\mathbb{Z}}_2$-linear combinations  of the $3$ generators of $I$. 
 Using Macaulay~2 \cite{M2}  we can verify that none of these  possible combinations produce a reduction of $I$. Therefore, $I$ is an $\m$-primary ideal with no $2$-generated reductions. Moreover, using Macaulay~2 again we find the integral closure of $I$ to be $\overline{I}=(I, x^3+xy^2)$ and thus $I$ is not integrally closed.}
\end{example}

Next we use Example~\ref{m-primary no red} to construct an integrally closed, $\m$-primary ideal in a regular local ring $(R,\m)$  that does not have $2$-generated reductions. This answers the question of Heinzer and Shannon \cite[Question~5.6]{HeS} discussed in Remark~\ref{Question HeS} in the negative.

\begin{example}\label{m-primary int closed no red}
{\rm Let $R=\mathbb{Z}_2[x,y]_{(x,y)}$ be a regular local ring with $\m=(x,y)$  the unique maximal ideal of $R$. 
Let $I$ be as in Example~\ref{m-primary no red}.  Recall that $\overline{I}=(I, x^3+xy^2)$. Since $I$ is ${\m}$-primary, so is $\overline{I}$.  We claim that   $\overline{I}$ is  does not have $2$-generated reductions. As in Example~\ref{m-primary no red}, we may assume any reduction of $\overline{I}$ is generated by  ${\mathbb{Z}}_2$-linear combinations of 
 the $4$ generators of $\overline{I}$. Using Macaulay~2 \cite{M2}  we  verify that none of these linear combinations produce a 2-generated reduction of $\overline{I}$. Therefore, $\overline{I}$ is an integrally closed $\m$-primary  ideal with no $2$-generated reductions.
}
\end{example}

\section*{acknowledgements}
We thank Bill Heinzer and Bernd Ulrich for suggesting Examples~\ref{Bill}~and~\ref{Bernd}, which served as our original motivation for this article.

\end{document}